\documentclass[11pt,a4paper]{amsart}
\usepackage{amssymb,amsrefs,mathrsfs,lineno}
\usepackage[matrix,arrow]{xy}
\usepackage{hyperref}

\usepackage{mathtools}

\renewcommand\aa{{\mathfrak a}}
\newcommand\ee{{\mathfrak e}}
\renewcommand\gg{{\mathfrak g}}
\renewcommand\qq{{\mathfrak q}}
\renewcommand\AA{{\mathfrak A}}
\newcommand\EE{{\mathfrak E}}
\newcommand\GG{{\mathfrak G}}
\newcommand\QQ{{\mathfrak Q}}
\def\U{{\mathbb U}}

\newcommand\Coalg{\mathbf{Coalg}_{\mathrm{coc}}}

\newcommand\Vect{\mathbf{Vect}}

\newtheorem{maintheorem}{Theorem}
\newtheorem{maincorollary}[maintheorem]{Corollary}
\newtheorem{theorem}{Theorem}[section]
\newtheorem{proposition}[theorem]{Proposition}
\newtheorem{corollary}[theorem]{Corollary}
\newtheorem{lemma}[theorem]{Lemma}

\theoremstyle{remark}

\theoremstyle{definition}

\def\presuper#1#2%
  {\mathop{}%
   \mathopen{\vphantom{#2}}^{#1}%
   \kern-\scriptspace%
   #2}

\begin{document}
\title{Wreath products of cocommutative Hopf algebras}
\author{Laurent Bartholdi}
\address{Mathematisches Institut, Georg-August Universit\"at, G\"ottingen, Germany}
\email{laurent.bartholdi@gmail.com}
\date{2 September 2026}
\begin{abstract}
  We define wreath products of cocommutative Hopf algebras, and show
  that they enjoy a universal property for coalgebra-split extensions,
  analogous to the Kaloujnine-Krasner theorem for groups.

  We show that the group ring of a wreath product of groups is the
  wreath product of their group rings.  The
  enveloping algebra of the wreath product of two Lie algebras is the
  irreducible component containing the identity in the wreath product
  of their enveloping algebras.

  We recover the aforementioned result that group extensions may be
  classified as certain subgroups of a wreath product, and that Lie
  algebra extensions may also be classified as certain subalgebras of
  a wreath product.
\end{abstract}
\maketitle


\section{Introduction}
Let $A,Q$ be cocommutative Hopf algebras. We construct the
\emph{wreath product} $A\wr Q$ of $A$ and $Q$, and show that it
satisfies a universal property with respect to containing extensions of $A$ by $Q$. The definition is very simple, in terms of
\emph{coalgebra exponentials}, see~\S\ref{ss:measuring}:
\[A\wr Q := A^Q\#Q.
\]
The notation $A^Q$ denotes the exponential in the cartesian closed
category of cocommutative coalgebras.  Thus it comes with a coalgebra
map $@\colon A^Q\otimes Q\to A$ and is characterized by
\begin{equation}\label{eq:coalg}
 \Coalg(C,A^Q)\cong\Coalg(C\otimes Q,A).
\end{equation}
We write $x@q$ only as shorthand for this evaluation map; in general
an element of $A^Q$ is not determined by its evaluations at ordinary
elements $q\in Q$.  The group-object structure of $A$ induces a
cocommutative Hopf algebra structure on $A^Q$, and $Q$ acts on it by
right translation.  These constructions, and the resulting smash
product, are detailed in Sections~\ref{ss:measuring}
and~\ref{ss:smash}.

Our first main result is that the wreath product of Hopf algebras
classifies their coalgebra-split extensions:
\begin{maintheorem}[Generalized Kaloujnine-Krasner theorem]\label{thm:KK}
  Let $A,Q$ be cocommutative Hopf algebras over a field. Denote by
  $\tau\colon A\wr Q\to Q$ the natural projection. Then the map
  $E\mapsto E$ defines a bijection between
  \[\frac{\left\{E:\xymatrix{{\Bbbk}\ar[r] & {A}\ar[r] & {E}\ar[r] & {Q}\ar[r] & {\Bbbk}}\text{ coalgebra-split}\right\}}{\text{equivalence of extensions}}\]
  and
  \[\frac{\left\{\begin{array}{l@{}r}&E\leq A\wr Q\text{ Hopf subalgebras}:\tau|_E\text{ has a coalgebra section},\\
                                     &E^{\operatorname{co}Q}=E\cap A^Q\cong A\text{ qua Hopf algebras via evaluation at }1
  \end{array}\right\}}{\text{conjugation by normalized coalgebra maps $Q\to A$.}}\]
  A normalized coalgebra map $\phi\colon Q\to A$ transposes
via~\eqref{eq:coalg} with $C=\Bbbk$ to a group-like
element $\widehat\phi\in A^Q$, and by \emph{conjugation by $\phi$} we mean the
ordinary inner Hopf automorphism of $A\wr Q$ given by conjugation by
$\widehat\phi\#1$.
\end{maintheorem}

Extensions of groups --- and of Hopf algebras --- with \emph{abelian}
kernel and a fixed action of $Q$ on $A$ are classified by the
cohomology group $H^2(Q,A)$;
see~\citelist{\cite{singer:extensions}\cite{agore-bontea-militaru:extensions}}. Kaloujnine
and Krasner considered wreath products as a means to classify
arbitrary extensions.

There are two fundamental examples of cocommutative Hopf algebras: the group ring
$\Bbbk\GG$ of a group $\GG$, with coproduct $\Delta(g)=g\otimes g$ for
all $g\in\GG$; and the universal enveloping algebra $\U(\gg)$ of a Lie
algebra $\gg$, with coproduct $\Delta(x)=x\otimes1+1\otimes x$ for all
$x\in\gg$.

Wreath products of groups have been considered since the beginnings
of group theory~\cite{jordan:subs}*{\S II.I.41}. The wreath product
$\AA\wr\QQ$ may be defined as the semidirect product
$\AA^\QQ\rtimes\QQ$; its universal property of containing all
extensions of $\AA$ by $\QQ$ is known as the Kaloujnine-Krasner
theorem. We show that the group ring of $\AA\wr\QQ$ is the wreath
product of the group rings of $\AA$ and $\QQ$, recovering in this
manner the Kaloujnine-Krasner theorem:

\begin{maintheorem}[Group rings]\label{thm:G}
  If $A=\Bbbk\AA$ and $Q=\Bbbk\QQ$ are group rings, then $A\wr Q\cong
  \Bbbk(\AA\wr\QQ)$ qua Hopf algebras.
\end{maintheorem}

\begin{maincorollary}[Kaloujnine-Krasner]\label{cor:Gkk}
  Let $\AA,\QQ$ be groups. Denote by $\tau\colon
  \AA\wr \QQ\to \QQ$ the natural projection. Then the map $\EE\mapsto \EE$
  defines a bijection between
  \[\frac{\left\{\EE:\xymatrix{{1}\ar[r] & {\AA}\ar[r] & {\EE}\ar[r] & {\QQ}\ar[r] & {1}}\right\}}{\text{equivalence of extensions}}\]
  and
  \[\frac{\left\{\begin{array}{l@{}r}&\EE\le \AA\wr \QQ\text{ subgroups}:\tau(\EE)=\QQ\text{ and }\EE\cap \AA^\QQ\cong \AA\\
        & \text{ via the evaluation map }\AA^\QQ\to \AA,f\mapsto f(1)\end{array}\right\}}{\text{conjugation by $(f,1)$, where $f\in\AA^\QQ$ and $f(1)=1$.}}\]
\end{maincorollary}

Special cases of wreath products of Lie algebras were considered in
various places in the
literature~\cites{jurman:glamc3,caranti-m-n:glamc1,netreba-s:wreath,dipietro:phd,shmelkin:wreath,sullivan:wreath,bondarenko:wreath}.
For the comparison with enveloping algebras we assume that $\Bbbk$ has
characteristic zero; in positive characteristic, the Lie algebras and universal enveloping algebras should be restricted.

The wreath product $\aa\wr\qq$ may be defined as
$\Vect(\U(\qq),\aa)\rtimes\qq$. An analogue of the Kaloujnine-Krasner
theorem was proved in~\cite{petrogradsky-r-s:wreath}. We show that the
universal enveloping algebra of $\aa\wr\qq$ is the identity component
of corresponding Hopf wreath product, recovering in this manner
the Kaloujnine-Krasner theorem:

\begin{maintheorem}[Lie algebras]\label{thm:L}
  Consider $A=\U(\aa)$ and $Q=\U(\qq)$, and let $(A\wr Q)_1$ denote the irreducible component
  of $A\wr Q$ containing its identity.  Then
  \[
      (A\wr Q)_1\cong\U(\aa\wr\qq)
  \]
  qua Hopf algebras.  In general $(A\wr Q)_1$ is a proper Hopf
  subalgebra of $A\wr Q$.
\end{maintheorem}

Here conjugation of Lie subalgebras means the restriction to primitive
elements of the Hopf conjugations defined above.  Every such Hopf
automorphism fixes the unit and therefore preserves the identity
component.  This is the precise sense in which the quotient in
Corollary~\ref{cor:Lkk} is taken.

\begin{maincorollary}[Kaloujnine-Krasner for Lie algebras, see~\cite{petrogradsky-r-s:wreath}]\label{cor:Lkk}
  Let $\aa,\qq$ be Lie algebras over a field.
  Denote by $\tau\colon
  \aa\wr \qq\to \qq$ the natural projection. Then the map $\ee\mapsto \ee$
  defines a bijection between
  \[\frac{\left\{\ee:\xymatrix{{0}\ar[r] & {\aa}\ar[r] & {\ee}\ar[r] & {\qq}\ar[r] & {0}}\right\}}{\text{equivalence of extensions}}\]
  and
  \[\frac{\left\{\begin{array}{l@{}r}&\ee\le \aa\wr \qq\text{ subalgebras}:\tau(\ee)=\qq\text{ and }\ee\cap \Vect(\U\qq,\aa)\cong \aa\\
        & \text{ via the evaluation map }\Vect(\U\qq,\aa)\to \aa,f\mapsto f(1)\end{array}\right\}}{\text{the conjugations described above.}}\]
\end{maincorollary}

\subsection{Assumptions}
Throughout Sections~\ref{ss:measuring}--\ref{sec:groups}, $\Bbbk$ is a
field and all coalgebras and Hopf algebras are cocommutative.
Section~\ref{sec:lie}
additionally assumes $\operatorname{char}(\Bbbk)=0$.  We write
$\Delta(u)=\sum u_1\otimes u_2$ in Sweedler notation, suppressing the
summation sign when no ambiguity can result.  Similarly, expressions
involving $u_1,u_2,u_3$ come from applying an iterated coproduct to
$u$.  We use Greek letters for morphisms and Latin letters for
elements.  Our convention for composition is
$(\phi\psi)(x)=\phi(\psi(x))$.

We use~\cite{sweedler:ha} and~\cite{milnor-moore:hopf} as general
references for Hopf algebras, and~\cite{montgomery:haaar} for
extensions of Hopf algebras.

\subsection{Thanks}
The author is grateful to Olivier Siegenthaler and Todd Trimble for
having contributed significant parts of this article --- he regrets
very much that they do not wish to be coauthors, since he values their
contributions as large as his. He is also grateful to the referee for
a very helpful and informative report that helped improve the text
very much.

\subsection{AI usage}
The original text was prepared in a pre-AI era, with no AI input. ChatGPT was
used to revise substantial parts of the text, following some referee
remarks. All mathematical content was reviewed, checked, and
substantially revised by the author, who is responsible for the final
text.

\section{Coalgebra exponentials}\label{ss:measuring}
Let $\Coalg$ denote the category of cocommutative counital coalgebras
over $\Bbbk$.  Its categorical product is the tensor product: the
projections $C\otimes D\to C,D$ are $1\otimes\varepsilon$ and
$\varepsilon\otimes1$, and the pairing of $\rho\colon R\to C$ and
$\sigma\colon R\to D$ is $(\rho\otimes\sigma)\Delta_R$.  Barr proved that this
category is complete and cartesian closed
in~\cite{barr:coalgebras}*{Theorem~5.3}.  Consequently, for
$C,D\in\Coalg$ there is an exponential coalgebra $D^C$ and a natural
bijection
\begin{equation}\label{eq:exponential}
 \Coalg(R,D^C)\xrightarrow{\sim}\Coalg(R\otimes C,D),\qquad
 \phi\longmapsto @(\phi\otimes1_C).
\end{equation}
The map $@\colon D^C\otimes C\to D$ is the counit of this adjunction;
in particular it is a coalgebra map.  We write $u@c$ for its value on
$u\otimes c$.

Formula~\eqref{eq:exponential}, rather than ordinary pointwise
evaluation, is the meaning of the exponential.  Here is how we use it:
two coalgebra maps $R\to D^C$ are equal if and only if their transposes
$R\otimes C\to D$ are equal.  By contrast, the linear map
$D^C\to\Vect(C,D)$ obtained by evaluating ordinary elements need not
be injective.

Since $@$ is a coalgebra map, its compatibility with counits and
coproducts reads explicitly
\begin{align}
 \varepsilon_D(u@c)&=\varepsilon_{D^C}(u)\varepsilon_C(c),
 \label{eq:eval-counit}\\
 \Delta_D(u@c)&=\sum (u_1@c_1)\otimes(u_2@c_2).
 \label{eq:eval-coproduct}
\end{align}
\subsection{Hopf algebra structure}
Let now $A$ be a cocommutative Hopf algebra and $C$ a cocommutative
coalgebra.  A cocommutative Hopf algebra is exactly a group object in
the cartesian category $\Coalg$.  Exponentiation of a group object is
again a group object, so $A^C$ is canonically a cocommutative Hopf
algebra.  Its structure maps are the unique coalgebra maps whose
transposes give the pointwise formulas
\begin{align}
 (xy)@c&=\sum(x@c_1)(y@c_2), & 1@c&=\varepsilon(c)1,
 \label{eq:exp-product}\\
 (Sx)@c&=S(x@c).&&\label{eq:exp-antipode}
\end{align}
All maps used here are coalgebra maps: cocommutativity makes
$\Delta_C\colon C\to C\otimes C$ a coalgebra map and makes the
antipode of $A$ a coalgebra map.  This is exactly where the
cocommutativity assumptions enter.

If $Q$ is a cocommutative Hopf algebra, it acts on $A^Q$ by right
translation.  The action $\star\colon Q\otimes A^Q\to A^Q$ is the
coalgebra map defined, via~\eqref{eq:exponential}, by
\begin{equation}\label{eq:translation}
       (q\star x)@r=x@(rq).
\end{equation}
Associativity of $Q$ proves that this is an action, while
\eqref{eq:exp-product} and the bialgebra identity for $Q$ give
\begin{align*}
 q\star1&=\varepsilon(q)1,&
 q\star(xy)&=\sum(q_1\star x)(q_2\star y).
\end{align*}
Hence $A^Q$ is a left $Q$-module Hopf algebra.

\section{Extensions of Hopf algebras}
For a Hopf map $\pi\colon E\to Q$, put
\[
 E^{\operatorname{co}Q}
 =\{e\in E\mid e_1\otimes\pi(e_2)=e\otimes1\}.
\]
An \emph{exact sequence} of cocommutative Hopf algebras is a diagram
\begin{equation}\label{eq:extension}
 \xymatrix{{\Bbbk}\ar[r]&A\ar[r]^-\iota&E\ar[r]^-\pi&Q\ar[r]&\Bbbk}
\end{equation}
in which $\iota$ is injective, $\pi$ is surjective,
$\iota(A)=E^{\operatorname{co}Q}$, and
$\ker\pi=E\iota(A)^+$.  See for example
\cites{andruskiewitsch-devoto:extensions,andruskiewitsch:extensions}
for this standard notion.  Note that a cocommutative Hopf algebra over
a field is faithfully flat over each Hopf subalgebra by
\cite{takeuchi:correspondence}*{Theorem~3.1}.

An equivalence of two sequences with the same $A$ and $Q$ is, by
definition, a Hopf algebra \emph{isomorphism} $E\to E'$ inducing the
identity on $A$ and $Q$.  Equivalently, one may require only a Hopf
algebra morphism $E\to E'$ inducing the identity on $A$ and $Q$: it
is automatically an isomorphism by the Short Five Lemma in the
semi-abelian category of cocommutative Hopf algebras over a
field~\cite{gran-sterck-vercruysse:semiabelian}*{Theorem~2.10}.

The sequence~\eqref{eq:extension} is \emph{coalgebra-split} if $\pi$
has a coalgebra section $\gamma\colon Q\to E$.  We call a normalized
one, with $\gamma(1)=1$, a \emph{cleavage}.  Any coalgebra section can
be normalized by left multiplication by $\gamma(1)^{-1}\in\iota(A)$.
Moreover it is automatically right $Q$-colinear, since
\[
 \gamma(q)_1\otimes\pi(\gamma(q)_2)=\gamma(q_1)\otimes q_2,
\]
and convolution-invertible, with
$\gamma^\dagger=S_E\gamma$.  Hence a coalgebra-split extension is
cleft in the usual sense of
\cite{montgomery:haaar}*{Section~7.2}.

The normal-basis isomorphism associated with a cleavage is
\begin{equation}\label{eq:normal-basis}
 \theta_\gamma\colon A\otimes Q\xrightarrow{\sim}E,
 \qquad a\otimes q\longmapsto\iota(a)\gamma(q),
\end{equation}
with inverse
\begin{equation}\label{eq:normal-basis-inverse}
 \theta_\gamma^{-1}(e)
 =\iota^{-1}\bigl(e_1\gamma^\dagger(\pi(e_2))\bigr)\otimes\pi(e_3).
\end{equation}
The expression in parentheses is coinvariant; multiplication of the
two displayed maps and the identities
$\gamma^\dagger*\gamma=\gamma*\gamma^\dagger=\eta\varepsilon$ prove
that they are inverse.  This is also the normal-basis part of the
cleft-extension theorem
\cite{montgomery:haaar}*{Theorem~8.2.4}.

\subsection{Smash and wreath products}\label{ss:smash}
Let $Q$ and $H$ be cocommutative Hopf algebras and let
$\star\colon Q\otimes H\to H$ be a coalgebra map satisfying
\begin{alignat}{2}
 q\star1&=\varepsilon(q)1,&
 q\star(hk)&=(q_1\star h)(q_2\star k),\label{eq:hopfaction:1}\\
 1\star h&=h,&q\star(r\star h)&=(qr)\star h.
 \label{eq:hopfaction:2}
\end{alignat}
The smash product $H\#Q$ has underlying coalgebra $H\otimes Q$ and
\begin{align}
 (h\#q)(k\#r)&=h(q_1\star k)\#q_2r,
 \label{eq:smash-product}\\
 S(h\#q)&=(S(q_1)\star S(h))\#S(q_2).
 \label{eq:smash-antipode}
\end{align}
The action being a coalgebra map gives
$\Delta(q\star h)=q_1\star h_1\otimes q_2\star h_2$.  Comparing
\begin{align*}
 \Delta((h\#q)(k\#r))
 &=h_1(q_1\star k_1)\#q_3r_1
   \otimes h_2(q_2\star k_2)\#q_4r_2,\\
 \Delta(h\#q)\Delta(k\#r)
 &=h_1(q_1\star k_1)\#q_2r_1
   \otimes h_2(q_3\star k_2)\#q_4r_2
\end{align*}
shows that cocommutativity of $Q$ makes the two expressions equal.
The Hopf axioms and~\eqref{eq:smash-antipode} then follow directly;
see also~\cite{molnar:sdp}*{Section~2}.

Apply this to the translation action~\eqref{eq:translation}.  We
define
\[
              A\wr Q=A^Q\#Q,
\]
and write $\tau(x\#q)=\varepsilon(x)q$.  This gives the split exact
sequence
\[
 \xymatrix{{\Bbbk}\ar[r]&A^Q\ar[r]&A\wr Q\ar[r]^-\tau&Q\ar[r]&\Bbbk},
\]
split by the Hopf map $q\mapsto1\#q$.

\subsection{Conjugation}
Let $\phi\colon Q\to A$ be a normalized coalgebra map.  Under
\eqref{eq:exponential}, with $R=\Bbbk$, it corresponds to a group-like
element $\widehat\phi\in A^Q$; we use a hat to denote this
correspondence.  Its inverse $S(\widehat\phi)$ has transpose
$\phi^\dagger=S_A\phi$.  We define
\[
 \operatorname{Ad}_\phi(z)=(\widehat\phi\#1)z(S\widehat\phi\#1).
\]
This is an ordinary inner Hopf automorphism and $\tau\operatorname{Ad}_\phi=\tau$.
For later calculations,~\eqref{eq:smash-product} gives
\begin{equation}\label{eq:conjugation}
 \operatorname{Ad}_\phi(x\#q)
 =\widehat\phi\,x\,(q_1\star S\widehat\phi)\#q_2.
\end{equation}
Applying $\operatorname{ev}_r\otimes1_Q$ to~\eqref{eq:conjugation} gives
\begin{equation}\label{eq:conjugation-evaluation}
 \bigl[\widehat\phi\,x\,(q_1\star S\widehat\phi)\bigr]@r\otimes q_2
 =\phi(r_1)(x@r_2)\phi^\dagger(r_3q_1)\otimes q_2.
\end{equation}

\section{The Kaloujnine-Krasner theorem for coalgebra-split extensions}
\label{sec:kkhopf}
We now prove Theorem~\ref{thm:KK}.

Call the Hopf subalgebras appearing on the right-hand side of
Theorem~\ref{thm:KK} admissible.

\begin{lemma}\label{lem:theta}
Let~\eqref{eq:extension} be coalgebra-split and let $\gamma$ be a
cleavage.  The formula
\begin{equation}\label{eq:theta-map}
 \theta(q,e)=\gamma(q_1)e_1
       \gamma^\dagger\bigl(q_2\pi(e_2)\bigr)
\end{equation}
defines a coalgebra map $\theta\colon Q\otimes E\to\iota(A)$.  Its
transpose $\widehat\theta\colon E\to A^Q$ is characterized by
\begin{equation}\label{eq:theta-transpose}
 \iota(\widehat\theta(e)@q)=\theta(q,e).
\end{equation}
It satisfies
\begin{align}
 \widehat\theta(ee')
 &=\widehat\theta(e_1)
   \bigl(\pi(e_2)\star\widehat\theta(e')\bigr),
 \label{eq:theta-product}\\
 \widehat\theta(1)&=1.\label{eq:theta-unit}
\end{align}
\end{lemma}
\begin{proof}
All operations in~\eqref{eq:theta-map} are coalgebra maps.  Here we use
cocommutativity to know that iterated diagonals and the antipode are
coalgebra maps.  A direct calculation of the right coaction gives
\begin{align*}
 \theta(q,e)_1\otimes\pi(\theta(q,e)_2)
 &=\gamma(q_1)e_1\gamma^\dagger(q_4\pi(e_4))\\
 &\qquad\otimes q_2\pi(e_2)S(q_3\pi(e_3))
 =\theta(q,e)\otimes1.
\end{align*}
Thus $\theta(q,e)\in E^{\operatorname{co}Q}=\iota(A)$, proving both
well-definedness and the coalgebra assertion.

By~\eqref{eq:exponential}, $\theta$ has the unique transpose
$\widehat\theta$.  To prove~\eqref{eq:theta-product}, transpose both
sides.  After applying~\eqref{eq:exp-product} and
\eqref{eq:translation}, the right-hand side becomes
\[
 \theta(q_1,e_1)\theta(q_2\pi(e_2),e').
\]
Inserting~\eqref{eq:theta-map}, the adjacent factors
$\gamma^\dagger(q_2\pi(e_2))\gamma(q_3\pi(e_3))$ cancel by
convolution, leaving $\theta(q,ee')$.  This proves the identity by the
uniqueness in~\eqref{eq:exponential}; the unit identity is similar.
\end{proof}

Define
\begin{equation}\label{eq:alpha}
 \alpha_\gamma\colon E\longrightarrow A\wr Q,
 \qquad
 \alpha_\gamma(e)=\widehat\theta(e_1)\#\pi(e_2).
\end{equation}
It is a coalgebra map, since it is the composite of $\Delta_E$ with
$\widehat\theta\otimes\pi$.  Equations~\eqref{eq:theta-product},
\eqref{eq:theta-unit}, and~\eqref{eq:smash-product} show that it is an
algebra map.  A bialgebra map between Hopf algebras automatically
commutes with antipodes, because both composites are convolution
inverses of the same map.  Hence $\alpha_\gamma$ is a Hopf map and
$\tau\alpha_\gamma=\pi$.

It is also injective.  Indeed,
\begin{equation}\label{eq:alpha-left-inverse}
 (\operatorname{ev}_1\otimes1)\alpha_\gamma(e)
 =\iota^{-1}\bigl(e_1\gamma^\dagger(\pi(e_2))\bigr)\otimes\pi(e_3)
 =\theta_\gamma^{-1}(e)
\end{equation}
by~\eqref{eq:normal-basis-inverse}.  If $a\in A$, then
$\alpha_\gamma(\iota(a))\in A^Q$ and
$\operatorname{ev}_1\alpha_\gamma(\iota(a))=a$.  It follows that
\[
 \alpha_\gamma(E)^{\operatorname{co}Q}
 =\alpha_\gamma(E)\cap A^Q=\alpha_\gamma(\iota(A)),
\]
and evaluation at $1$ identifies this Hopf algebra with $A$.

\begin{lemma}\label{lem:change-cleavage}
The conjugacy class of $\alpha_\gamma(E)$ is independent of the
cleavage and of the representative of the extension.
\end{lemma}
\begin{proof}
For two normalized cleavages $\gamma,\gamma'$ define
$\phi\colon Q\to A$ by
\[
 \iota(\phi(q))=\gamma'(q_1)\gamma^\dagger(q_2).
\]
The same coaction computation as in Lemma~\ref{lem:theta} shows that
the right-hand side is coinvariant.  Cocommutativity shows that $\phi$ is a
normalized coalgebra map, and
$\gamma'(q)=\iota(\phi(q_1))\gamma(q_2)$.  Transposing the definitions
and using~\eqref{eq:conjugation-evaluation} gives
\[
       \alpha_{\gamma'}=\operatorname{Ad}_\phi\alpha_\gamma.
\]
If $\chi\colon E\to E'$ is an equivalence of extensions, then
$\chi\gamma$ is a cleavage of $E'$ and the construction gives
$\alpha_{\chi\gamma}\chi=\alpha_\gamma$.  Comparing $\chi\gamma$
with any chosen cleavage of $E'$ proves the assertion.
\end{proof}

We have therefore constructed, from every coalgebra-split extension,
an admissible Hopf subalgebra of $A\wr Q$, well defined up to the
conjugations specified in Theorem~\ref{thm:KK}.

Conversely, let $\lambda\colon B\hookrightarrow A\wr Q$ be an admissible
Hopf subalgebra and put $\pi=\tau\lambda$.  By hypothesis $\pi$ has a
coalgebra section and
$K=B^{\operatorname{co}Q}=B\cap A^Q$ is identified with $A$ by
$\operatorname{ev}_1$.  Takeuchi's correspondence gives
$\ker\pi=BK^+$; hence we have a coalgebra-split exact sequence
\[
 \xymatrix{{\Bbbk}\ar[r]&A\ar[r]&B\ar[r]^-\pi&Q\ar[r]&\Bbbk}.
\]
Conjugation by a normalized $\phi$ restricts to an equivalence of these
extensions because it fixes $\tau$ and evaluation at $1$ on $K$.

We verify both inverse identities.  Starting with
an extension $E$, the first construction produces
$\alpha_\gamma(E)$ and the inverse construction returns the same
extension through the isomorphism $\alpha_\gamma$; hence this
composition is the identity on equivalence classes.

For the other composition, define the coalgebra and right-comodule map
\[
 \mu=(\operatorname{ev}_1\otimes1)\lambda\colon B\to A\otimes Q.
\]
Choose any normalized cleavage $\gamma$ of $B$ and define
\[
 \rho=(1\otimes\varepsilon)\mu\gamma\colon Q\to A.
\]
The map $\mu$ is left $A$-linear and right $Q$-colinear, so
\eqref{eq:normal-basis} gives
\[
 \mu(\iota(a)\gamma(q))=a\rho(q_1)\otimes q_2
\]
and $\rho$ is a normalized coalgebra map.  The map $\mu$ is invertible,
with inverse in these coordinates
$a\otimes q\mapsto\iota(a\rho^\dagger(q_1))\gamma(q_2)$.  Thus
\[
 \gamma'(q)=\mu^{-1}(1\otimes q)
\]
is a normalized cleavage.  Let
$\beta=(1\otimes\varepsilon)\lambda\colon B\to A^Q$ be the base component of
$\lambda$.  The comodule identity gives
\begin{equation}\label{eq:base-reconstruction}
       \lambda(b)=\beta(b_1)\#\pi(b_2).
\end{equation}
Using $\mu\gamma'(q)=1\otimes q$ and the smash-product multiplication, one
computes
\[
 \operatorname{ev}_1\!\left(
 \lambda(\gamma'(r_1))\lambda(b_1)\lambda(\gamma'{}^\dagger(r_2\pi(b_2)))
 \right)=\beta(b)@r.
\]
The left-hand side is precisely the transpose~\eqref{eq:theta-transpose} defining
the map $\widehat\theta$ associated with the cleavage $\gamma'$.  Hence
$\widehat\theta=\beta$ by~\eqref{eq:exponential}, and~\eqref{eq:alpha} and
\eqref{eq:base-reconstruction} give $\alpha_{\gamma'}=\lambda$.  Any other
cleavage gives a conjugate image by Lemma~\ref{lem:change-cleavage}.
This proves the second inverse identity and Theorem~\ref{thm:KK}.

\section{Groups}\label{sec:groups}
We recall the universal property of
wreath products of groups mentioned in the introduction:
\begin{theorem}[Kaloujnine-Krasner, \cite{kaloujnine-krasner:extensions}]\label{thm:KKg}
  Let $\EE$ be an extension of $\AA$ by $\QQ$:
  \[\xymatrix{{1}\ar[r] & {\AA}\ar[r] & {\EE}\ar[r]^{\pi} & {\QQ}\ar[r] & {1}}.\]
  Then there is an injective homomorphism
  $\alpha\colon\EE\to\AA\wr\QQ$ such that $\tau\alpha=\pi$.

  Conversely, if $\EE$ is a subgroup of $\AA\wr\QQ$ which maps onto $\QQ$ by
  the natural map $\tau\colon \AA\wr\QQ\to\QQ$, and such that $\ker\tau\cap\EE$
  is isomorphic to $\AA$ via $f\mapsto f(1)$, then $\EE$ is an extension
  of $\AA$ by $\QQ$.
\end{theorem}

Although the proof is classical, we cannot resist including it, since
it is particularly short, and is essentially the proof of
Theorem~\ref{thm:KKg}:
\begin{proof}[Sketch of proof]
  Choose a set-theoretic section $q\mapsto\widetilde q$ of $\pi$.
  We define $\alpha\colon \EE\to\AA\wr\QQ$ by
  \[\alpha(e)=\left(q\mapsto \widetilde q e(\widetilde{q\pi(e)})^{-1},\pi(e)\right).\]
  It is clear that $\alpha$ is injective, and an easy check shows that
  $\alpha$ is a homomorphism. Conversely, if $\EE$ is a subgroup of
  $\AA\wr\QQ$ as in the statement of the theorem, then $\pi=\tau|_\EE$
  defines the extension.
\end{proof}

\subsection{Proof of Theorem~\ref{thm:G}}
The wreath product of groups $\AA,\QQ$ is the semidirect product
$\AA^\QQ\rtimes\QQ$, and the group ring of a semidirect product is a
smash product of group rings.  It therefore suffices to identify the
coalgebra exponential.
\begin{proposition}
  Let $X,Y$ be sets, and let $\Bbbk X,\Bbbk Y$ be their group-like
  coalgebras, with $\Delta(x)=x\otimes x$ and $\varepsilon(x)=1$ for
  all $x\in X$; and similarly for $Y$.

  Then the coalgebras $(\Bbbk Y)^{\Bbbk X}$ and $\Bbbk(Y^X)$ are
  isomorphic.
\end{proposition}

\begin{proof}
  It is enough to compare the functors on finite-dimensional
  cocommutative coalgebras, since every coalgebra is the filtered union
  of its finite-dimensional subcoalgebras and coalgebra maps are
  determined by their restrictions to these subcoalgebras.  Write such
  a coalgebra as $R^*$, where $R$ is a finite-dimensional commutative
  algebra.  A coalgebra map $R^*\to\Bbbk Y$ is equivalent
  to a finite-support family $(e_y)_{y\in Y}$ of pairwise orthogonal
  idempotents of $R$ with sum $1$.  Call such a family a
  $Y$-distribution.

  By~\eqref{eq:exponential}, a map
  $R^*\to(\Bbbk Y)^{\Bbbk X}$ is an $X$-tuple
  $(e_{x,y})_{y\in Y}$ of $Y$-distributions.  All these idempotents
  commute because $R$ is commutative.  The Boolean algebra of
  idempotents of a finite-dimensional algebra is finite, so the
  partitions $(e_{x,y})_{y\in Y}$ have a finite common refinement by
  nonzero pairwise orthogonal idempotents $p_1,\dots,p_n$ with sum $1$.
  Each $p_i$ selects, for every $x$, a unique $y$ and hence a function
  $\phi_i\colon X\to Y$.  Combining equal $\phi_i$ gives a
  finite-support $Y^X$-distribution.

  Conversely, a $Y^X$-distribution $(e_\phi)$ has marginals
  \[
       e_{x,y}=\sum_{\phi(x)=y}e_\phi.
  \]
  The two constructions are inverse and natural in $R$.  Hence the
  two coalgebras represent the same functor and are isomorphic.
\end{proof}

Under this isomorphism, the Hopf operations and the translation action
are the pointwise group operations and right translation.  Taking
$X=\QQ$ and $Y=\AA$, then forming the smash product, proves
Theorem~\ref{thm:G}.

\subsection{Proof of Corollary~\ref{cor:Gkk}}

By $\mathscr G(A)$ we denote the \emph{group-like} elements of a Hopf
algebra $A$, defined as
\[\mathscr G(A)=\{x\in A\colon \Delta(x)=x\otimes x\text{ and }\varepsilon(x)=1\}.\]
\begin{lemma}
  Let $A$ be a Hopf algebra. Then $\mathscr G(A)$ is linearly
  independent in $A$. The following are equivalent:
  \begin{enumerate}
  \item $A$ is a group algebra;
  \item the canonical Hopf map $\Bbbk\mathscr G(A)\to A$ is an
    isomorphism;
  \item $\mathscr G(A)$ is a linear basis of $A$.
  \end{enumerate}
\end{lemma}
\begin{proof}
  We argue by induction.  Suppose $x_1,\dots,x_n$ are independent
  group-like elements and a group-like $x$ lies in their span, say
  $x=\sum_i c_ix_i$. Then
  \[\sum_i c_i x_i\otimes x_i = \Delta(x) = x\otimes x = \sum_{i,j}c_ic_j x_i\otimes x_j.\]
  Therefore $c_ic_j=0$ for all $i\neq j$, and $c_i^2=c_i$ for all $i$,
  so, over a field and using $\varepsilon(x)=1$, exactly one $c_i$ is
  $1$.  Thus $x$ is one of the $x_i$.  A minimal linear dependence
  among group-like elements would contradict this observation, proving
  linear independence.  The equivalence follows.
\end{proof}

\begin{corollary}\label{cor:gprings}
  Let $A,Q$ be the group rings of groups $\AA,\QQ$ respectively. Then
  there is a bijection between equivalence classes of coalgebra-split
  extensions of $A$ by $Q$ and equivalence classes of group extensions
  of $\AA$ by $\QQ$, which relates each extension of $\AA$ by $\QQ$ to
  its group ring.
\end{corollary}
\begin{proof}
  Consider first an extension
  \[\xymatrix{{1}\ar[r] & {\AA}\ar[r]^{\iota} & {\EE}\ar[r]^{\pi} & {\QQ}\ar[r] & {1}},
  \]
  and set $E=\Bbbk\EE$. Then the natural maps $\Bbbk\iota\colon A\to E$
  and $\Bbbk\pi\colon E\to Q$ turn $E$ into an extension of $A$ by $Q$,
  which is coalgebra-split because any set-theoretic section of $\pi$
  extends linearly to a coalgebra section.

  Conversely, consider a coalgebra-split extension
  \begin{equation}\label{eq:gprings:1}
    \xymatrix{{\Bbbk}\ar[r] & {A}\ar[r]^{\iota} & {E}\ar[r]^{\pi} &
      {Q}\ar[r] & {\Bbbk}},
  \end{equation}
  and choose a cleavage $\gamma\colon\Bbbk\QQ\to E$.  Every
  $\gamma(q)$ is group-like.  By the normal-basis
  isomorphism~\eqref{eq:normal-basis}, the elements
  $\iota(a)\gamma(q)$, for $a\in\AA$ and $q\in\QQ$, form a basis of
  $E$ consisting of group-like elements.  Thus
  $E=\Bbbk\mathscr G(E)$.  Set $\EE=\mathscr G(E)$. Then the restriction
  $\overline\iota\colon \AA\to\EE$ is injective because $\iota$ is
  injective, and the restriction $\overline\pi\colon \EE\to\QQ$ is
  surjective because $\pi$ is split qua coalgebra map. We certainly
  have $\overline\pi\,\overline\iota(a)=1$ for every $a\in\AA$,
  because~\eqref{eq:gprings:1} is exact. Finally, consider
  $e\in\ker\overline\pi$; then
  $e\in E^{\operatorname{co}Q}\cap\EE=\iota(\AA)$, so
  \[\xymatrix{{1}\ar[r] & {\AA}\ar[r]^{\overline\iota} &
    {\EE}\ar[r]^{\overline\pi} & {\QQ}\ar[r] & {1}}
  \]
  is exact and its group algebra is the original Hopf extension.
\end{proof}

\noindent Corollary~\ref{cor:Gkk} now follows from
Theorems~\ref{thm:KK} and~\ref{thm:G}, and
Corollary~\ref{cor:gprings}.

\section{Lie algebras}\label{sec:lie}
In this section, we assume for simplicity that $\Bbbk$ has characteristic zero; if $\Bbbk$ has positive characteristic, the Lie algebras should be assumed to be restricted, and the universal enveloping algebra should also be restricted.  Let $\aa$ and $\qq$
be Lie algebras and put
\[
  \aa^\qq=\Vect(\U(\qq),\aa).
\]
The coproduct of $\U(\qq)$ makes $\aa^\qq$ a Lie algebra, with
\[
  [f,g](u)=\sum [f(u_1),g(u_2)].
\]
There is a left action of $\qq$ on $\aa^\qq$ by derivations,
\[
  (q\star f)(u)=f(uq).
\]
Indeed, the identity
$\Delta(uq)=\sum u_1q\otimes u_2+\sum u_1\otimes u_2q$ gives
$q\star[f,g]=[q\star f,g]+[f,q\star g]$.  The Lie wreath product is
\[
  \aa\wr\qq=\aa^\qq\rtimes\qq.
\]
Thus, if its elements are written as $f\oplus q$, then
\begin{equation}\label{eq:wrbracket}
 [f\oplus q,g\oplus r]
 =\left(u\mapsto\sum[f(u_1),g(u_2)]+g(uq)-f(ur)\right)
   \oplus[q,r].
\end{equation}
For a cocommutative Hopf algebra $H$, write
\[
 \mathscr P(H)=\{x\in H:\Delta(x)=x\otimes1+1\otimes x\}.
\]
Also let $H_1$ denote the irreducible (equivalently, connected)
component containing $1$.  It is the largest connected Hopf
subalgebra of $H$.

\begin{proposition}\label{prop:primitive-exponential}
  For cocommutative Hopf algebras $A,Q$, evaluation induces a natural
  Lie-algebra isomorphism
  \[
       \mathscr P(A^Q)\xrightarrow{\ \sim\ }
       \Vect(Q,\mathscr P(A)),\qquad
       x\longmapsto (q\mapsto x@q).
  \]
\end{proposition}
\begin{proof}
  Let $D=\Bbbk1\oplus\Bbbk d$ be the pointed coalgebra with
  $\Delta(d)=d\otimes1+1\otimes d$.  A primitive element of $A^Q$ is
  the same thing as a coalgebra map $D\to A^Q$ which sends $1$ to the
  unit of $A^Q$.  By the exponential adjunction, such a map is the
  same thing as a coalgebra map $D\otimes Q\to A$ whose restriction to
  $1\otimes Q$ is $q\mapsto\varepsilon(q)1$.  Its restriction to
  $d\otimes Q$ has the form $d\otimes q\mapsto f(q)$, and the
  coalgebra identities say precisely that
  \[
       \Delta f(q)=f(q)\otimes1+1\otimes f(q),\qquad
       \varepsilon f(q)=0.
  \]
  Hence $f$ is an arbitrary linear map $Q\to\mathscr P(A)$.  The
  pointwise multiplication formula in $A^Q$ shows that the induced
  bracket is
  $[f,g](q)=\sum[f(q_1),g(q_2)]$.
\end{proof}

\begin{proof}[Proof of Theorem~\ref{thm:L}]
  Apply Proposition~\ref{prop:primitive-exponential} to
  $A=\U(\aa)$ and $Q=\U(\qq)$.  Since
  $\mathscr P(\U(\aa))=\aa$, it gives
  \[
       \mathscr P(A^Q)=\Vect(\U(\qq),\aa)=\aa^\qq.
  \]
  The Cartier-Kostant-Milnor-Moore theorem for connected
  cocommutative Hopf algebras
  (see~\cite{milnor-moore:hopf}*{Theorem~5.18}) therefore gives
  \[
       (A^Q)_1\cong\U(\aa^\qq).
  \]
  The smash product has underlying coalgebra $A^Q\otimes Q$ and $Q$ is
  connected, so
  \[
       (A\wr Q)_1=(A^Q)_1\#Q.
  \]
  Under the preceding primitive-element identification, the
  right-translation action of $Q$ differentiates to
  $(q\star f)(u)=f(uq)$.  Finally, the standard enveloping-algebra
  isomorphism
  \[
       \U(\aa^\qq\rtimes\qq)
       \cong \U(\aa^\qq)\#\U(\qq)
  \]
  proves the assertion.

  The qualification by the identity component is necessary:
  group-like elements of $A^Q$ correspond to coalgebra maps
  $Q\to A$, and these need not consist only of the unit map.
\end{proof}

\begin{corollary}\label{cor:lierings}
  Let $A=\U(\aa)$ and $Q=\U(\qq)$.  Passage to universal enveloping
  algebras induces a bijection between equivalence classes of Lie
  algebra extensions of $\aa$ by $\qq$ and equivalence classes of
  coalgebra-split Hopf algebra extensions of $A$ by $Q$.
\end{corollary}
\begin{proof}
  Given a Lie extension
  \[
    \xymatrix{{0}\ar[r] & {\aa}\ar[r] & {\ee}\ar[r] & {\qq}\ar[r] & {0}},
  \]
  the PBW theorem shows that applying $\U$ gives an exact Hopf
  extension.  A linear section $\qq\to\ee$, together with ordered PBW
  bases, extends to a normalized coalgebra section
  $\U(\qq)\to\U(\ee)$, so this extension is coalgebra-split.

  Conversely, let
  \[
    \xymatrix{{\Bbbk}\ar[r] & {A}\ar[r]^{\iota} & {E}\ar[r]^{\pi} &
      {Q}\ar[r] & {\Bbbk}}
  \]
  be coalgebra-split.  The normal-basis
  isomorphism~\eqref{eq:normal-basis} identifies the coalgebra $E$
  with $A\otimes Q$.  Hence $E$ is connected, and the
  Cartier-Kostant-Milnor-Moore theorem identifies it canonically
  with $\U(\mathscr P(E))$.  A normalized coalgebra section sends
  primitives of $Q$ to primitives of $E$, so
  $\mathscr P(E)\to\qq$ is surjective.  Moreover, if
  $x\in\mathscr P(E)$ and $\pi(x)=0$, then
  $(1\otimes\pi)\Delta(x)=x\otimes1$; exactness therefore puts
  $x$ in $\iota(A)$, and primitiveness then puts it in $\iota(\aa)$.
  Thus
  \[
     0\longrightarrow\aa\longrightarrow\mathscr P(E)
       \longrightarrow\qq\longrightarrow0
  \]
  is exact, and its enveloping-algebra extension is the given one.
  These constructions preserve and reflect equivalences and are
  mutually inverse.
\end{proof}

\begin{theorem}[Petrogradsky-Razmyslov-Shishkin]\label{thm:kklie}
  Every extension of $\aa$ by $\qq$ embeds in $\aa\wr\qq$ over $\qq$.
  Conversely, if $\ee\leq\aa\wr\qq$, the natural projection
  $\tau_\ee\colon\ee\to\qq$ is surjective, and evaluation at $1$ restricts
  to an isomorphism
  \[
       \ker(\tau_\ee)\xrightarrow{\ \sim\ }\aa,
  \]
  then $\ee$ is an extension of $\aa$ by $\qq$.
\end{theorem}
\begin{proof}
  The first assertion follows from Theorem~\ref{thm:KK},
  Corollary~\ref{cor:lierings}, and Theorem~\ref{thm:L}, by taking
  primitive elements in the identity component.  The converse is
  immediate from the stated kernel and surjectivity conditions.  This
  is the Lie-algebra Kaloujnine-Krasner theorem
  of~\cite[Theorem~3.1]{petrogradsky-r-s:wreath}.
\end{proof}

\noindent Corollary~\ref{cor:Lkk} follows from
Theorem~\ref{thm:KK}, Corollary~\ref{cor:lierings}, and the
primitive-element description in
Proposition~\ref{prop:primitive-exponential}.



\begin{bibdiv}
\begin{biblist}
\bib{agore-bontea-militaru:extensions}{article}{
  author={Agore, Ana-Loredana},
  author={Bontea, Costel-Gabriel},
  author={Militaru, Gigel},
  title={Classifying coalgebra split extensions of Hopf algebras},
  journal={J. Algebra Appl.},
  volume={12},
  date={2013},
  number={5},
  pages={1250227, 24},
  issn={0219-4988},
  review={\MR {3055582}},
  doi={10.1142/S0219498812502271},
}

\bib{andruskiewitsch:extensions}{article}{
  author={Andruskiewitsch, Nicol\'as},
  title={Notes on extensions of Hopf algebras},
  journal={Canad. J. Math.},
  volume={48},
  date={1996},
  number={1},
  pages={3\ndash 42},
  issn={0008-414X},
  review={\MR {1382474 (97c:16046)}},
  doi={10.4153/CJM-1996-001-8},
}

\bib{andruskiewitsch-devoto:extensions}{article}{
  author={Andruskiewitsch, Nicol\'as},
  author={Devoto, Jorge A.},
  title={Extensions of Hopf algebras},
  journal={Algebra i Analiz},
  volume={7},
  date={1995},
  number={1},
  pages={22\ndash 61},
  issn={0234-0852},
  translation={ journal={St. Petersburg Math. J.}, volume={7}, date={1996}, number={1}, pages={17\ndash 52}, issn={1061-0022}, },
  review={\MR {1334152 (96f:16044)}},
}

\bib{barr:coalgebras}{article}{
  author={Barr, Michael},
  title={Coalgebras over a commutative ring},
  journal={J. Algebra},
  volume={32},
  date={1974},
  number={3},
  pages={600\ndash 610},
  doi={10.1016/0021-8693(74)90161-6},
}

\bib{bondarenko:wreath}{article}{
  author={Bondarenko, Natalia~V.},
  title={Lie algebras associated with wreath products of elementary abelian groups},
  language={English, with English and Russian summaries},
  journal={Mat. Stud.},
  volume={26},
  date={2006},
  number={1},
  pages={3\ndash 16},
  issn={1027-4634},
  review={\MR {2285063 (2007m:17026)}},
}

\bib{caranti-m-n:glamc1}{article}{
  author={Caranti, Andrea},
  author={Mattarei, Sandro},
  author={Newman, Michael~F.},
  title={Graded Lie algebras of maximal class},
  journal={Trans. Amer. Math. Soc.},
  volume={349},
  date={1997},
  number={10},
  pages={4021\ndash 4051},
  issn={0002-9947},
  review={\MR {1443190 (98a:17027)}},
  doi={10.1090/S0002-9947-97-02005-9},
}

\bib{gran-sterck-vercruysse:semiabelian}{article}{
  author={Gran, Marino},
  author={Sterck, Florence},
  author={Vercruysse, Joost},
  title={A semi-abelian extension of a theorem by Takeuchi},
  journal={J. Pure Appl. Algebra},
  volume={223},
  date={2019},
  number={10},
  pages={4171--4190},
  issn={0022-4049},
  review={\MR {3958087}},
  doi={10.1016/j.jpaa.2019.01.004},
}

\bib{jordan:subs}{book}{
  author={Jordan, Camille},
  title={Trait\'e des substitutions et des \'equations alg\'ebriques},
  language={French},
  series={Les Grands Classiques Gauthier-Villars. [Gauthier-Villars Great Classics]},
  note={Reprint of the 1870 original},
  publisher={\'Editions Jacques Gabay},
  place={Sceaux},
  date={1989},
  pages={xvi+670},
  isbn={2-87647-021-7},
  review={\MR {1188877 (94c:01039)}},
}

\bib{jurman:glamc3}{article}{
  author={Jurman, Giuseppe},
  title={Graded Lie algebras of maximal class. III},
  journal={J. Algebra},
  volume={284},
  date={2005},
  number={2},
  pages={435\ndash 461},
  issn={0021-8693},
  review={\MR {2114564 (2005k:17041)}},
  doi={10.1016/j.jalgebra.2004.11.006},
}

\bib{kaloujnine-krasner:extensions}{article}{
  author={Kaloujnine, Lev~A.},
  author={Krasner, Marc},
  title={Le produit complet des groupes de permutations et le probl\`eme d'extension des groupes},
  language={French},
  journal={C. R. Acad. Sci. Paris},
  volume={227},
  date={1948},
  pages={806\ndash 808},
  review={\MR {0027758 (10,351e)}},
}

\bib{milnor-moore:hopf}{article}{
  author={Milnor, John~W.},
  author={Moore, John~C.},
  title={On the structure of Hopf algebras},
  journal={Ann. of Math. (2)},
  volume={81},
  date={1965},
  pages={211\ndash 264},
  issn={0003-486X},
  review={\MR {0174052 (30 \#4259)}},
}

\bib{molnar:sdp}{article}{
  author={Molnar, Richard~K.},
  title={Semi-direct products of Hopf algebras},
  journal={J. Algebra},
  volume={47},
  date={1977},
  number={1},
  pages={29\ndash 51},
  issn={0021-8693},
  review={\MR {0498612 (58 \#16700)}},
}

\bib{montgomery:haaar}{book}{
  author={Montgomery, Susan},
  title={Hopf algebras and their actions on rings},
  series={CBMS Regional Conference Series in Mathematics},
  volume={82},
  publisher={Published for the Conference Board of the Mathematical Sciences, Washington, DC},
  date={1993},
  pages={xiv+238},
  isbn={0-8218-0738-2},
  review={\MR {1243637 (94i:16019)}},
}

\bib{netreba-s:wreath}{article}{
  author={Su{\v s}{\v c}ans{\cprime }ki{\u \i }, Vitali{\u \i }~{\=I}},
  author={Netreba, Nataliya~V.},
  title={Wreath product of Lie algebras and Lie algebras associated with Sylow $p$-subgroups of finite symmetric groups},
  journal={Algebra Discrete Math.},
  date={2005},
  number={1},
  pages={122\ndash 132},
  issn={1726-3255},
  review={\MR {2148825}},
}

\bib{petrogradsky-r-s:wreath}{article}{
  author={Petrogradski{\u \i }, Victor~M.},
  author={Razmyslov, Yuri P.},
  author={Shishkin, E. O.},
  title={Wreath products and Kaluzhnin-Krasner embedding for Lie algebras},
  journal={Proc. Amer. Math. Soc.},
  volume={135},
  date={2007},
  number={3},
  pages={625\ndash 636 (electronic)},
  issn={0002-9939},
  review={\MR {2262857 (2008a:17020)}},
  doi={10.1090/S0002-9939-06-08502-9},
}

\bib{dipietro:phd}{thesis}{
  author={Di Pietro, Cristina},
  title={Wreath products and modular Lie algebras},
  date={2005},
  organization={Universit\`a degli Studi di L'Aquila},
  pages={ix+89},
}

\bib{singer:extensions}{article}{
  author={Singer, William~M.},
  title={Extension theory for connected Hopf algebras},
  journal={J. Algebra},
  volume={21},
  date={1972},
  pages={1\ndash 16},
  issn={0021-8693},
  review={\MR {0320056 (47 \#8597)}},
}

\bib{shmelkin:wreath}{article}{
  author={{\v {S}}mel{\cprime }kin, Alfred~L.},
  title={Wreath products of Lie algebras, and their application in group theory},
  language={Russian},
  note={Collection of articles commemorating Aleksandr Gennadievi\v c Kuro\v s},
  journal={Trudy Moskov. Mat. Ob\v s\v c.},
  volume={29},
  date={1973},
  pages={247\ndash 260},
  issn={0134-8663},
  review={\MR {0379612 (52 \#517)}},
}

\bib{sullivan:wreath}{article}{
  author={Sullivan, Frances~E.},
  title={Wreath products of Lie algebras},
  journal={J. Pure Appl. Algebra},
  volume={35},
  date={1985},
  number={1},
  pages={95\ndash 104},
  issn={0022-4049},
  review={\MR {772163 (86h:17002)}},
}

\bib{sweedler:ha}{book}{
  author={Sweedler, Moss~E.},
  title={Hopf algebras},
  series={Mathematics Lecture Note Series},
  publisher={W. A. Benjamin, Inc., New York},
  date={1969},
  pages={vii+336},
  review={\MR {0252485 (40 \#5705)}},
}

\bib{takeuchi:correspondence}{article}{
  author={Takeuchi, Mitsuhiro},
  title={A correspondence between Hopf ideals and sub-Hopf algebras},
  journal={Manuscripta Math.},
  volume={7},
  date={1972},
  number={3},
  pages={251--270},
  doi={10.1007/BF01579722},
}

\end{biblist}
\end{bibdiv}
\end{document}